\documentclass[a4paper]{amsart}

\usepackage{amscd}
\usepackage{amssymb}
\usepackage{amsxtra}
\usepackage{dbnsymb}

\allowdisplaybreaks[2]

\theoremstyle{plain}
\newtheorem{thm}{Theorem}[section]
\newtheorem{cor}[thm]{Corollary}
\newtheorem{prop}{Proposition}
\newtheorem{lem}[thm]{Lemma}

\theoremstyle{definition}
\newtheorem{defn}{Definition}

\theoremstyle{remark}
\newtheorem*{rmk}{Remark}
\newtheorem*{rmks}{Remarks}
\newtheorem*{exmp}{Example}
\newtheorem*{exmps}{Examples}

\renewcommand{\epsilon}{\varepsilon}
\renewcommand{\kappa}{\varkappa}
\renewcommand{\theta}{\vartheta}

\newcommand{\IC}{\ensuremath{\mathbb C}}

\newcommand{\IN}{\ensuremath{\mathbb N}}
\newcommand{\IQ}{\ensuremath{\mathbb Q}}

\newcommand{\IZ}{\ensuremath{\mathbb Z}}

\newcommand{\ch}{\ensuremath{\mathrm{ch}}}
\newcommand{\td}{\ensuremath{\mathrm{td}}}

\newcommand{\Alt}{\ensuremath{\mathrm{Alt}}}

\newcommand{\End}{\ensuremath{\mathrm{End}}}

\newcommand{\Sym}{\ensuremath{\mathrm{Sym}}}
\newcommand{\Gr}{\ensuremath{\mathrm{Gr}}}
\DeclareMathOperator{\Pic}{Pic}

\newcommand{\card}[1]{\ensuremath{\left|#1\right|}}
\newcommand{\varcard}[1]{\ensuremath{\operatorname{card}}}
\newcommand{\ord}[1]{\ensuremath{\operatorname{ord}}}
\newcommand{\set}[1]{\ensuremath{\left\{#1\right\}}}

\newcommand{\tr}{\ensuremath{\operatorname{tr}}}

\newcommand{\HH}{\ensuremath{\mathrm H}}
\newcommand{\KK}{\ensuremath{\mathrm K}}

\DeclareMathOperator{\pr}{pr}
\DeclareMathOperator{\id}{id}

\newcommand{\shO}{\ensuremath{\mathcal O}}
\newcommand{\OX}{\ensuremath{\mathcal O_X}}
\newcommand{\Tang}{\ensuremath{\mathcal T}}
\newcommand{\TangX}{\ensuremath{\Tang_X}}

\newcommand{\RW}{\ensuremath{\operatorname{RW}}}

\def\trivalent{\mathrm{t}}

\newcommand{\BB}{\ensuremath{\mathcal B}}
\newcommand{\hBB}{\ensuremath{\hat{\BB}}}
\newcommand{\BBo}{\ensuremath{\BB^\trivalent}}
\newcommand{\hBBo}{\ensuremath{\hBB^\trivalent}}
\newcommand{\hl}{\ensuremath{\hat\ell}}

\newcommand{\SA}[1]{\ensuremath{A^{({#1})}}}
\newcommand{\KA}[1]{\ensuremath{K^{#1}\!A}}
\newcommand{\HA}[1]{\ensuremath{A^{[{#1}]}}}

\newcommand{\mf}[1]{\mathfrak{#1}}

\newcommand{\To}{\:\longrightarrow\:}

\numberwithin{equation}{section}

\def\N{\IN}
\def\Z{\IZ}
\def\Q{\IQ}

\def\C{\IC}
\def\eps{\epsilon}

\begin{document}

\title[HRR on an irreducible symplectic K\"ahler]
{Hirzebruch-Riemann-Roch formulae on irreducible symplectic K\"ahler
  manifolds}
\author[M.\ Britze]{Michael Britze}
\author[M.\ A.\ Nieper]{Marc A.\ Nieper}
\address{Mathematisches Institut der Univ.\ zu K\"oln \\
        Weyertal 86--90 \\ 50931 K\"oln \\ Germany}
\email{mbritze@mi.uni-koeln.de \\ mail@marc-nieper.de}
\thanks{We are very grateful to Daniel Huybrechts for having carefully
  read preliminary versions of this paper and to Daniel Huybrechts,
  Manfred Lehn and many others for their support to us and helpful
  discussions about the subject.}

\begin{abstract}
  In this article we investigate Hirzebruch-Riemann-Roch formulae for
  line bundles on irreducible symplectic K\"ahler manifolds. As Huybrechts
  has shown in ~\cite{huybrechts99}, for every irreducible
  complex K\"ahler manifold $X$ of dimension $2n$, there are
  numbers $a_0, a_2, \dots, a_{2n}$ such that $\chi(L) = \sum_{k =
    0}^n a_{2k}/(2k)!  q_X(c_1(L))^k$ for the Euler characteristic of
  a line bundle $L$, where $q_X: \HH^2(X, \C) \to \C$ is the
  Beauville-Bogomolov quadratic form of $X$.

  Using Rozansky-Witten classes similar to Hitchin and Sawon
  in~\cite{hitchin99}, we obtain a formula expressing the $a_{2k}$ in
  terms of Chern numbers of $X$.  Furthermore, for the $n$-th
  generalized Kummer variety $\KA n$ (see~\cite{beauville83}), we
  prove $\chi(L) = (n + 1) \binom{q(c_1(L)) / 2 + n} n$ by purely
  algebro-geometric methods, where $q$ is the form $q_X$ up to a
  positive rational constant. A similar formula is already known for
  the Hilbert scheme of zero-dimensional subschemes of length $n$ on a
  K3-surface (c.f.~\cite{lehn99}).

  Using our results, we are able to calculate all Chern numbers of the
  generalized Kummer varieties $\KA n$ for $n \leq 5$. For $n \leq 4$
  these results were previously obtained by Sawon
  (c.f.~\cite{sawon99}).
\end{abstract}

\maketitle

\tableofcontents

\section{Some linear algebra}

In this section, let $k$ denote a field of characteristic zero, $V$ a
$k$-vector space and $A$ a $k$-algebra.

We identify the exterior algebra $\bigwedge (V^*) \otimes_k A$ with
the space $\bigoplus_{r = 0}^n \Alt^r(V, A)$ of alternating
multilinear forms on $V$ with values in $A$ as vector spaces by setting
\begin{gather}
  \label{equ:formident}
  ((\alpha_1 \wedge \dots \wedge \alpha_r) \otimes a)(v_1, \dots, v_r)
  = \det((\alpha_i(v_j))_{ij}) \cdot a
\end{gather}
for $\alpha_1, \dots, \alpha_r \in V^*$, $v_1, \dots, v_r \in V$ and
$a \in A$.

\begin{defn}
  We call an element $\sigma \in \bigwedge^2(V^*)$ a \emph{symplectic
    form on $V$} and $V$ together with $\sigma$ a \emph{symplectic
    $k$-vector space} if the map
  \begin{gather}
    I_\sigma: V \to V^*, v \mapsto \sigma(\cdot, v)
  \end{gather}
  is an isomorphism of vector spaces. 
\end{defn}

\begin{rmk}
  A symplectic vector space is always of even dimension.
\end{rmk}

If $\sigma$ is a fixed symplectic form on $V$, we will identify $V$
and $V^*$ using the isomorphism $I_\sigma$. In particular, we have an
induced dual symplectic form $\tilde\sigma \in \bigwedge^2 V$ on $V^*$.

For the rest of this section, let $\sigma$ be a fixed symplectic form
on a $2n$-dimensional vector space $V$.

\begin{defn}
  We define a non-singular pairing
  \begin{gather}
    \left<\cdot, \cdot\right>: \left(\bigwedge V^* \otimes_k A\right)
    \otimes_k \left(\bigwedge V^* \otimes_k A\right) \to A
  \end{gather}
  by setting
  \begin{gather}
    \label{equ:formscp}
      \left<(\alpha_1 \wedge \dots \wedge \alpha_r) \otimes a,
        (\beta_1 \wedge \dots \wedge \beta_s) \otimes b\right>
      = \delta_{rs} \det ((\sigma^*(\alpha_i, \beta_j))_{ij}) \cdot ab
  \end{gather}
  for $\alpha_1, \dots, \alpha_r, \beta_1, \dots, \beta_s \in V^*$ and
  $a, b \in A$, where $\delta_{rs}$ is Kronecker's $\delta$.
\end{defn}

The following proposition will be used later in the case that $V$ is a
tangent space of a complex manifold at a point $p$ and $A$ is the algebra of
anti-holomorphic forms at $p$.

\begin{prop}
  \label{prop:laexp}
  Let $\sigma$ be a symplectic form on $V$ and $\int: \bigwedge(V^*)
  \otimes_k A \to \bigwedge^{2n}(V^*) \otimes_k A$ the canonical
  projection onto the forms of top degree.

  For every $\alpha \in \bigwedge (V^*)$ we have
  \begin{gather}
    \int(\alpha \wedge \exp\sigma) = \int(\left<\alpha,
      \exp\sigma\right>\exp\sigma).
  \end{gather}
\end{prop}

\begin{proof}
  We can assume that $A = k$ and $\alpha = \alpha_1 \wedge \dots
  \wedge \alpha_{2p} \in \bigwedge^{2p}(V^*)$ with $p \in \N_0$,
  $\alpha_i \in V^*$. So we have to prove
  \begin{equation*}
    \alpha^{2p} \wedge \sigma^{(n - p)}
    = \frac{(n - p)!}{p! \cdot n!} \left<\alpha, \sigma^p\right>\sigma^n.
  \end{equation*}
  Let $e_1, \dots, e_{2n}$ be a symplectic basis of $V$ and $\theta^1,
  \dots, \theta^{2n}$ the corresponding dual basis of $V^*$, i.e.\ 
  $\sigma = \sum_{i = 1}^n \theta^{2i - 1} \wedge \theta^{2i}$ and
  $\tilde\sigma = \sum_{i = 1}^n e_{2i - 1} \wedge e_{2i}$. It follows
  that $\sigma^n = n! \cdot \theta^1 \wedge \dots \wedge \theta^{2n}$.
  It is
  \begin{multline*}
    \left<\alpha \wedge \sigma^{(n - p)}, \sigma^n\right> \\
    \shoveleft{
      = n! (n - p)! \sum_{1 \leq i_1 < \dots < i_{n - p} \leq n}
      \begin{aligned}        
        & \bigl<\alpha_1 \wedge \dots \wedge \alpha_{2p} \wedge
        \theta^{2i_1 - 1} \wedge \theta^{2i_i} \wedge \dots \\
        & \indent \dots \wedge \theta^{2i_{n - p} - 1}
        \wedge \theta^{2i_{n - p}}, \theta^1 \wedge \dots \wedge
        \theta^{2n}
        \bigr>
      \end{aligned}
      }
    \\
    \shoveleft{
      = n! (n - p)! \sum_{1 \leq j_1 < \dots < j_p \leq n}
      \left<\alpha_1 \wedge \dots \wedge \alpha_{2p},
        \theta^{2j_1 - 1} \wedge \theta^{2j_1} \wedge \dots \wedge
        \theta^{2j_p - 1} \wedge \theta^{2j_p}
      \right>
      }
    \\
    \shoveleft{
      = \frac{n! (n - p)!}{p!} \left<\alpha, \sigma^p\right>.
      }
    \\ \ 
  \end{multline*}
  Since $\sigma^n$ spans $\bigwedge^{2n}(V^*)$ and
  $\left<\sigma^n, \sigma^n\right> = n!^2$, this proves the proposition.
\end{proof}

\begin{rmk}
  Hitchin and Sawon have stated this result for $\alpha$ being a
  two-form in~\cite{hitchin99}, p.\ 8. Note that they identify the
  exterior algebra over $V^*$ with the alternating forms on $V$ in a
  different way than we do.

  Extending their formula to arbitrary degree of $\alpha$ is crucial
  for this work.
\end{rmk}

\begin{defn}
  For every set $S$, we denote by $\bigwedge_k S$ the Gra\ss
  mann algebra generated by the elements of $S$ over $k$. If $S'$ is a
  subset of $S$, we view $\bigwedge_k(S')$ canonically as a
  subalgebra of $\bigwedge_k S$. We denote by $\Sym^n_k S$ the $n$-th
  symmetric product of the $k$-vector space spanned by the elements of
  $S$.
\end{defn}


\section{Graph homology}

\subsection{The graph homology space}

In this article, \emph{graph} means a collection of vertices connected
by edges, i.e.\ every edge connects two vertices.  We want to call a
half-edge (i.e.\ an edge together with an adjacent vertex) of a graph
a \emph{flag}. So, every edge consists of exactly two flags. Every
flag belongs to exactly one vertex of the graph. On the other hand, a
vertex is given by the set of its flags. It is called
\emph{univalent}, if there is only one flag belonging to it, and it is
called \emph{trivalent}, if there are exactly three flags belonging to
it. We shall identify edges and vertices with the set of their flags.

A graph is called \emph{vertex-oriented} if, for every vertex, a
cyclic ordering of its flags is fixed.

\begin{defn}
  A \emph{Jacobi diagram} is a vertex-oriented graph with only uni-
  and trivalent vertices. A \emph{connected Jacobi diagram} is a
  Jacobi diagram which is connected as a graph. A \emph{trivalent
  Jacobi diagram} is a Jacobi diagram with no univalent vertices.
\end{defn}

We define the \emph{degree of a Jacobi diagram} to be the number of
its vertices. It is always an even number.

\begin{exmp}
  The empty graph is a Jacobi diagram, denoted by $1$.
  The unique Jacobi diagram consisting of two univalent vertices (which
  are connected by an edge) is denoted by $\ell$.
\end{exmp}

We can always draw a Jacobi diagram in a planar drawing so that it
looks like a planar graph with vertices of valence $1$, $3$ or
$4$. Each $4$-valent vertex has to be interpreted as a crossing of two
non-connected edges of the drawn graph and not as one of its vertices.
Further, we want the counter-clockwise ordering of the flags at each
trivalent vertex in the drawing to be the same as the given cyclic
ordering.

In drawn Jacobi diagrams, we also use a notation like $\cdots\overset n
-\cdots$ for a part of a graph which looks like a long line with $n$
univalent vertices (``legs'') attached to it, for example
$\ldots\bot\!\bot\!\bot\ldots$ for $n = 3$. The position of
$n$ indicates the placement of the legs relative to the ``long line''.

\begin{defn}
  Let $S$ be a totally ordered set with $n$ elements. We set
  \begin{gather}
    \epsilon(S) := s_1 \wedge \dots \wedge s_n \in \bigwedge_\Q\nolimits S.
  \end{gather}
  Here, $s_1, \dots, s_n$ are the elements of $S$ in increasing order.
\end{defn}

\begin{defn}
  Let $S$ be a cyclicly ordered set with an odd number $n$ of
  elements. We set
  \begin{gather}
    \epsilon(S) := s_1 \wedge \dots \wedge s_n \in \bigwedge_\Q\nolimits S.
  \end{gather}
  Here, $s_1, \dots, s_n$ are the elements of $S$ in an order
  compatible with the given cyclic one. The definition of
  $\epsilon(S)$ does not depend on this order as long as the
  compatibility condition is fulfilled.
\end{defn}

\begin{exmp}
  If $t$ is a trivalent vertex in a Jacobi diagram, it makes sense to
  write $\epsilon(t)$ because we have said that we identify a vertex
  with the set of the flags belonging to it.
\end{exmp}

\begin{defn}
  Let $\Gamma$ be a Jacobi diagram with $k$ trivalent and $l$
  univalent vertices, so $m := \frac{3k + l} 2$ is the number of
  its edges. Let $F$ be the set of its flags, $E$ the set of
  its edges, $T$ the set of its trivalent vertices and $U$ the set of
  its univalent vertices.
  
  A choice of total orderings of the sets $T$, $U$ and every set $e
  \in E$ (recall that an edge is identified with the set of the two
  flags belonging to it) are said to be \emph{compatible with the
    orientation of $\Gamma$} if the equality
  \begin{gather}
    \epsilon(t_1) \wedge \dots \wedge \epsilon(t_k)
    \wedge \epsilon(u_1) \wedge \dots \wedge \epsilon(u_l)
    = \epsilon(e_1) \wedge \dots \wedge \epsilon(e_m)
  \end{gather}
  holds in the Gra\ss mann algebra $\bigwedge_\Q F$ generated by the
  elements of $F$.
  Here, $t_1, \dots, t_k$ are the elements of $T$ in increasing order,
  $u_1, \dots, u_l$ are the elements of $U$ in increasing order, and
  $e_1, \dots, e_m$ are all the edges in arbitrary order.
\end{defn}

\begin{defn}
  We define $\BB$ to be the $\Q$-vector space spanned by all Jacobi
  diagrams modulo the IHX relation
  \begin{gather}
    \IGraph = \HGraph - \XGraph
  \end{gather}
  and the anti-symmetry (AS) relation
  \begin{gather}
    \YGraph + \TwistedY = 0,
  \end{gather}
  which can be applied anywhere within a diagram. (For this definition
  see also~\cite{barnatan95} and ~\cite{thurston00}.)

  Furthermore, let $\BB'$ be the subspace of $\BB$ spanned by all Jacobi
  diagrams not containing $\ell$ as a component, and let $\BBo$ be the
  subspace of $\BB'$ spanned by all trivalent Jacobi diagrams. All
  these are graded and double-graded. The grading is induced by the
  degree of the Jacobi diagrams, the double-grading by the number of
  univalent and trivalent vertices.

  The completion of $\BB$ (resp.\ $\BB'$, resp.\ $\BBo$) with respect
  to the grading will be denoted by $\hBB$ (resp.\ $\hBB'$, resp.\ $\hBBo$).
\end{defn}

All these spaces are called \emph{graph homology spaces} and their
elements are called \emph{graph homology classes} or \emph{graphs} for
short. 

\begin{rmks}
  The subspaces $\BB_k$ of $\hBB$ spanned by the Jacobi diagrams of
  degree $k$ are always of finite dimension. The subspace $\BB_0$ is
  one-dimensional and spanned by the graph homology class of the empty
  diagram $1$.
\end{rmks}

\begin{exmp}
  If $\gamma$ is a graph which has a part looking like
  $\cdots\overset n -\cdots$, it will become $(-1)^n \gamma$ if we
  substitute the part $\cdots\overset n -\cdots$ by $\cdots\underset n
  -\cdots$ due to the anti-symmetry relation.
\end{exmp}

\subsection{Operations with graphs and special graphs}

\begin{defn}
  Disjoint union of Jacobi diagrams induces a bilinear map
  \begin{gather}
    \hBB \times \hBB \to \hBB, (\gamma, \gamma') \mapsto \gamma \cup
    \gamma'.
  \end{gather}
\end{defn}

By mapping $1 \in \Q$ to $1 \in \hBB$, the space $\hBB$ becomes a
graded $\Q$-algebra, which has no components in odd degrees. Often, we
omit the product sign ``$\cup$''. $\BB$, $\BB'$, $\BBo$ and so on
are subalgebras.

\begin{defn}
  Let $k \in \N$. We call the graph homology class of the Jacobi
  diagram $\overset{2k}{\bigcirc}$ the \emph{$2k$-wheel $w_{2k}$},
  i.e.\ $w_2 = \twowheel$, $w_4 = \fourwheel$ and so on. It has $2k$
  univalent and $2k$ trivalent vertices. Further, we set $w_0 := 1$.

  For $k_1, k_2 \in \N$, we call the graph homology class of the
  Jacobi diagram $\overset{k_1}{\underset{k_2}{\ThetaGraph}}$ a 
  \emph{double-wheel}, denoted by $w_{k_1, k_2}$. In particular,
  $w_{0, 0} = \ThetaGraph$.
\end{defn}

Next, we define a very special element $\Omega$ of the graph homology
space with very remarkable properties, which have been proven by
Bar-Natan, Le and Thurston (see~\cite{thurston00}).

\begin{defn}
  Let the series $(b_{2k})_k \in \C^{\N_0}$ of the \emph{modified
    Bernoulli numbers} (c.f.~\cite{thurston00}) be defined by
  \begin{equation}
    \sum_{k = 0}^\infty b_{2k} x^{2k} = \frac 1 2 \ln\frac{\sinh
      \frac x 2}{\frac x 2}.
  \end{equation}
  Let $\Omega \in \hBB'$ be the image of the element
  \begin{equation}
    \exp\left(\sum_{k = 1}^\infty b_{2k} x_{2k}\right) \in
    \C[[(x_{2k})_{k \in \N}]]
  \end{equation}
  under the morphism $\C[[(x_{2k})]] \to \hBB' \otimes_\Q \C$ of
  $\C$-algebras that maps $x_{2k}$ to $w_{2k}$.

  For any $\mu \in \C$ we set $\Omega(\mu) := \sum_{k = 0}^\infty
  \Omega_k \mu^k$, where $\Omega_k$ is the homogeneous component of
  degree $2k$ of $\Omega$. (Note that $\Omega_k = 0$ for odd $k$.)
\end{defn}

The following remark is stated in~\cite{thurston00}.
\begin{rmk}
  The modified Bernoulli numbers are connected to the usual Bernoulli
  numbers $B_1, B_2, B_3, \dots$ via
  \begin{gather}
    b_{2k} = \frac{B_{2k}}{4k (2k)!}
  \end{gather}
  for all $k \in \N$. In addition to this, $b_0 = 0$.

  The generating function of the (usual) Bernoulli numbers is given by
  \begin{gather}
    \sum_{k = 0}^\infty \frac{B_k}{k!} t^k = \frac t{e^t - 1}.
  \end{gather}
  Note that $B_k = 0$ for $k > 1$ and $k$ odd. Furthermore, $B_0 = 1$ and
  $B_1 = -\frac 1 2$.
\end{rmk}

Let $\Gamma$ be a Jacobi diagram and $u,
u'$ be two different univalent vertices of $\Gamma$. At least one
of them should not belong to a component $\ell$ of $\Gamma$. Let $v$
(resp.\ $v'$) be the vertex $u$ (resp.\ $u'$) is attached to. The process of
\emph{glueing the vertices $u$ and $u'$} means to remove $u$ and $u'$
together with the edges connecting them to $v$ resp.\ $v'$ and to add
a new edge between $v$ and $v'$. Thus, we arrive at a new graph
$\Gamma/(u, u')$, whose number of trivalent vertices is the numbers of
trivalent vertices of $\Gamma$ and whose number of univalent vertices
is the numbers of univalent vertices of $\Gamma$ minus two. To make it
a Jacobi diagram we define the cyclic orientation of the flags at $v$
(resp.\ $v'$) to be the cyclic orientation of the flags at $v$ (resp.\ 
$v'$) in $\Gamma$ with the flag belonging to the edge connecting $v$
(resp.\ $v'$) with $u$ (resp.\ $u'$) replaced by the flag belonging to
the added edge.

For example, glueing the two univalent vertices of $w_2$ leads to the
graph $\ThetaGraph$.

Of course, the process of glueing to univalent vertices given above
does not work if $u$ and $u'$ are the two univalent vertices of
$\ell$, thus our assumption on $\Gamma$.

\begin{defn}
  Let $\Gamma, \Gamma'$ be two Jacobi diagrams, at least one of them
  without $\ell$ as a component and $U = \set{u_1, \dots, u_n}$
  resp.\ $U'$ the sets of their univalent vertices.  We define
  \begin{gather}
    \hat\Gamma(\Gamma') := \sum_{\substack{f: U \hookrightarrow U' \\
    \text{injective}}} (\Gamma \cup \Gamma')/(u_1, f(u_1))/\dots/(u_n,
    f(u_n)),
  \end{gather}
  viewed as an element in $\hBB$.

  This induces for every $\gamma \in \hBB$ a $\hBBo$-linear map
  \begin{gather}
    \hat\gamma: \hBB' \to \hBB', \gamma' \mapsto \hat\gamma(\gamma').
  \end{gather}
\end{defn}

\begin{exmp}
  Set $\partial := \frac 1 2 \hl$.
  It is is an endomorphismus of $\hBB'$ of degree $-2$. For example,
  $\partial \twowheel = \ThetaGraph$. By setting
  \begin{gather}
    \partial(\gamma, \gamma') := \partial(\gamma\cup\gamma') -
    \partial(\gamma)\cup\gamma' - \gamma\cup\partial(\gamma')
  \end{gather}
  for $\gamma, \gamma' \in \hBB'$,
  we have the following formula for all $\gamma \in
  \hBB'$:
  \begin{gather}
    \partial(\gamma^n) = \binom n 1 \partial(\gamma)\gamma^{n - 1} +
    \binom n 2 \partial(\gamma, \gamma) \gamma^{n - 2}.
  \end{gather}

  Acting by $\partial$ on a Jacobi diagram means to glue two of its
  univalent vertices in all possible ways, acting by $\partial(\cdot,
  \cdot)$ on two Jacobi diagrams means to connect them by glueing an
  univalent vertex of the first with an univalent vertex of the second
  in all possible ways. For example, we have
  \begin{gather}
    \partial(w_{2k}) = k \sum_{n = 0}^{2k - 2} w_{n, 2k - 2 - n}
  \end{gather}
  for $k \in \N$ and, by the IHX relation,
  \begin{gather}
    \partial(w_{2k_1}, w_{2k_2}) = 8 k_1 k_2 w_{2 k_1 - 1, 2 k_2 - 1} 
  \end{gather}
  for $k_1, k_2 \in \N$.
\end{exmp}

\begin{defn}
  Let $\Gamma, \Gamma'$ be two Jacobi diagrams, at least one of them
  without $\ell$ as a component, and $U = \set{u_1, \dots, u_n}$
  resp.\ $U'$ the sets of their univalent vertices.  We define
  \begin{gather}
    \left<\Gamma, \Gamma'\right>
    := \sum_{\substack{f: U \to U' \\
        \text{bijective}}} (\Gamma \cup \Gamma')/(u_1, f(u_1))/\dots/(u_n,
    f(u_n)),
  \end{gather}
  viewed as an element in $\hBBo$.

  This induces a $\hBBo$-bilinear map
  \begin{gather}
    \left<\cdot, \cdot\right>, \hBB' \times \hBB \to \hBBo,
  \end{gather}
  which is symmetric on $\hBB' \times \hBB'$.
\end{defn}

Note that $\left<\Gamma, \Gamma\right>$ is zero unless $\Gamma$ and
$\Gamma'$ have equal numbers of univalent vertices. In this case, the
expression is the sum over all possibilities to glue the univalent
vertices of $\Gamma$ with univalent vertices of $\Gamma'$.

\begin{prop}
  The map $\left<1, \cdot\right>: \hBB \to \hBBo$ is the canonical
  projection map, i.e.\ it removes all non-trivalent components from a
  graph. Furthermore, for $\gamma \in \hBB'$ and $\gamma' \in \hBB$, we
  have
  \begin{gather}
    \label{equ:ellandpartial}
    \left<\gamma, \frac\ell 2 \gamma'\right> = \left<\partial \gamma,
      \gamma'\right>.
  \end{gather}

  For $\gamma, \gamma' \in \hBB$, we have the following (combinatorial)
  formula:
  \begin{equation}
    \label{equ:scpandpartial}
    \left<\exp(\partial)(\gamma\gamma'), 1\right>
    = \left<\exp(\partial)\gamma, \exp(\partial)\gamma'\right>
  \end{equation}
\end{prop}

\begin{proof}
  The formula~\eqref{equ:ellandpartial} should be clear from the
  definitions. 

  Let us investigate~\eqref{equ:scpandpartial} a bit more.
  We can assume that $\gamma$ and $\gamma'$ are Jacobi diagrams with
  $l$ resp.\ $l'$ univalent vertices and $l + l' = 2n$ with $n \in \N_0$. So
  we have to prove
  \begin{equation*}
    \label{equ:combin}
    \frac{\partial^n}{n!} (\gamma\gamma')
    = \sum_{\substack{m, m' = 0 \\ l - m = l' - m'}}^\infty
      \left<\frac{\partial^m}{m!} \gamma, \frac{\partial^{m'}}{m'!}
        \gamma'\right>,
  \end{equation*}
  since $\left<\cdot, 1\right>: \hBB \to \hBBo$ means to remove the
  components with at least one univalent
  vertex. Recalling the meaning of $\left<\cdot, \cdot\right>$, it
  should be clear that~\eqref{equ:combin} follows from the
  fact that applying $\frac{\partial^k}{k!}$ on a Jacobi diagram means to
  glue $k$ pairs of its univalent vertices in all possible ways.
\end{proof}

\subsection{$\Omega$ as an eigenvector of the operator $\partial$}

As said above, the element $\Omega$ plays a central role in the
``Wheeling Theorem'' (see~\cite{thurston00}, where Bar-Natan, Le and
Thurston have gave a knot theoretical prove of this theorem). Hitchin
and Sawon (\cite{hitchin99}) discovered that this theorem together
with the ideas of Rozansky and Witten (\cite{rozansky97}) can be used
to deduce some interesting facts about characteristic classes on
irreducible symplectic K\"ahler manifolds.

To prove our Hirzebruch-Riemann-Roch formula on irreducible symplectic
K\"ahler manifolds along the ideas of Hitchin and Sawon, we could make
use of the ``Wheeling Theorem''. However, it tells us a lot
more about $\Omega$ than we actually need. In fact, we only need the
statement given in theorem~\ref{thm:omega} about $\Omega$ to prove our
results. Similarly, everything stated in~\cite{hitchin99} that is
based on the ``Wheeling Theorem'' can also
be based on theorem~\ref{thm:omega}.

The following theorem is a corollary of lemma~6.2
in~\cite{thurston00}. Nevertheless, we give another proof here, which
does not use any knot theory.

\begin{thm}
  \label{thm:omega}
  For each $\mu \in \C$, the graph $\Omega(\mu)$ is an eigenvector of
  the endomorphismus $\partial: \hBB' \to \hBB'$ with the eigenvalue
  $\frac {\mu^2}{48} \ThetaGraph$, i.e.\ 
  \begin{gather}
    \partial\Omega(\mu) = \frac{\mu^2}{48} \ThetaGraph \Omega(\mu).
  \end{gather}
\end{thm}

To make the proof of the theorem readable, we first state two
lemmata. The first lemma is a combinatorial one, which makes use of
the IHX relations. It collects some results of Chmutov, Dasbach and Duzhin
Dasbach (see~\cite{chmutov99} and \cite{dasbach98}).

\begin{lem}
  \label{lem:chmutov}
  Let $\mathcal W$ be the subspace of $\hBB$ that is spanned by all
  graphs $w_{i, j}$ with $i, j \in \N_0$.
  Let $P: \mathcal W \to \Sym^3_\Q((x_n)_{n \in \N_0})$ be the map
  defined by
  \begin{gather}
    \label{equ:chmutovp}
    P(w_{i, j}) =
    \begin{cases}
      2 \sum_{l, m = 0}^\infty (-1)^{l + m} \binom{i}{l}
      \binom{j}{m} x_l x_m x_{i + j - l - m} & \text{for $i + j$ even}
      \\
      0 & \text{for $i + j$ odd.}
    \end{cases}
  \end{gather}
  Then $P$ is injective.
\end{lem}

\begin{proof}
  First, we show that $P$ is well-defined.

  For every $n \in \N_0$, let $\mathcal W_n$ be the subspace of $\mathcal
  W$ spanned by all $w_{i, j}$ with $i + j = n$.
  It is enough to show that for all $n \in \N_0$ there exists a injective map
  $P_n: \mathcal W_n \to \Sym^3_\Q((x_n)_{n \in \N_0})$ that
  fulfills \eqref{equ:chmutovp}. Further, we can restrict ourselves to
  the case of even $n$ due to lemma 6.2 of~\cite{dasbach98}, which
  says that $w_{i, j}$ is homologous to zero for odd $n$. (This
  follows at once from the anti-symmetry relation.)
  
  For any $N > n$, the map $P_{\mathfrak{gl}(N)}$ defined in
  section~3.1 of~\cite{chmutov99} can be taken as $P$ when
  restricted to $\mathcal W_n$. This is because of proposition~4.5
  of~\cite{chmutov99}, where it is shown that $P_{\mathfrak{gl}(N)}$
  evaluated at $w_{i, j}$ equals the right hand side
  of~\eqref{equ:chmutovp}.

  It remains to show that $P_n$ is injective.
  This can be proven by a dimension argument: By lemma~6.2 and
  lemma~6.8 of~\cite{dasbach98} the image of $\mathcal W_n$ under
  $P_n$ has at least the dimension of $\mathcal W_n$,
  so $P_n$ is injective.
\end{proof}

\begin{lem}
  \label{lem:bernoulli}
  Let $B_0, B_1, B_2, \dots$ denote the Bernoulli numbers. The
  following formula holds in $\Sym^3_\Q((x_n)_{n \in \N_0})$:
  \begin{multline}
    \sum_{k = 2}^\infty \frac{B_k}{k!} \sum_{n = 0}^{k - 2}
    \sum_{l, m = 0}^\infty (-1)^{l + m} \binom n l \binom{k - 2 - n}m
    x_l x_m x_{k - 2 - l - m} \\
    + \sum_{i, j = 2}^\infty
    \frac{B_i}{i!}\frac{B_j}{j!} \sum_{l, m = 0}^\infty (-1)^{l + m}
    \binom{i - 1}l \binom{j - 1}m
    x_l x_m x_{i + j - 2 - l - m}
    = \frac 1 {12} x_0^3.
  \end{multline}
\end{lem}

\begin{proof}
  In $\Q[X_1, X_2, X_3]$, we calculate
  \begin{multline*}
    \shoveleft{
      \sum_{\pi \in \mathfrak S_3} \Biggl(
      \sum_{k = 2}^\infty \frac{B_k}{k!} \sum_{n = 0}^{k - 2}
      \sum_{l, m = 0}^\infty (-1)^{l + m} \binom n l \binom{k - 2 - n}m
      X_{\pi(1)}^l X_{\pi(2)}^m X_{\pi(3)}^{k - 2 - l - m}
      }
    \\
    \shoveleft{
      + \sum_{i, j = 2}^\infty
      \frac{B_i}{i!}\frac{B_j}{j!} \sum_{l, m = 0}^\infty (-1)^{l + m}
      \binom{i - 1}l \binom{j - 1}m
      X_{\pi(1)}^l X_{\pi(2)}^m X_{\pi(3)}^{i + j - 2 - l - m}
      \Biggr)
      }
    \\
    \shoveleft{
      = \sum_{\pi \in \mathfrak S_3} \Biggl(  
      \sum_{k = 2}^\infty \frac{B_k}{k!} \sum_{n = 0}^{k - 2}
      \left(X_{\pi(3)} - X_{\pi(1)}\right)^n
      \left(X_{\pi(3)} - X_{\pi(2)}\right)^{k - 2 - n}
      }
    \\
    \shoveleft{
      + \sum_{i, j = 2}^\infty
      \frac{B_i}{i!}\frac{B_j}{j!}
      (X_{\pi(3)} - X_{\pi(1)})^{i - 1} (X_{\pi(3)} - X_{\pi(2)})^{j -
        1} \Biggr)
      }
    \\
    \shoveleft{
      = \sum_{\pi \in \mathfrak S_3} \Biggl(
      \frac 1{X_{\pi(1)} - X_{\pi(2)}} \left(
        \begin{aligned}
          & \frac 1{X_{\pi(3)} - X_{\pi(2)}}
          \sum_{k = 2}^\infty \frac{B_k}{k!}
          \left(X_{\pi(3)} - X_{\pi(2)}\right)^k \\
          & \indent - \frac 1{X_{\pi(3)} - X_{\pi(1)}}
          \sum_{k = 2}^\infty \frac{B_k}{k!}
          \left(X_{\pi(3)} - X_{\pi(1)}\right)^k
        \end{aligned}
      \right)
      }
    \\
    \shoveleft{
      + \frac 1{(X_{\pi(3)} - X_{\pi(1)})(X_{\pi(3)} - X_{\pi(2)})}
      }
    \\
    \shoveleft{
      \cdot
      \left(
        \sum_{k = 2}^\infty \frac{B_k}{k!}
        \left(X_{\pi(3)} - X_{\pi(1)}\right)^k
      \right)
      \left(
        \sum_{k = 2}^\infty \frac{B_k}{k!}
        \left(X_{\pi(3)} - X_{\pi(2)}\right)^k
      \right) 
      \Biggr)
      }
    \\
    \shoveleft{
      = \sum_{\pi \in \mathfrak S_3} \Biggl(
      \frac 1{X_{\pi(1)} - X_{\pi(2)}}
      \left(
        \begin{aligned}
          & \frac 1{\exp\left(X_{\pi(3)} - X_{\pi(2)}\right) - 1}
          - \frac 1{X_{\pi(3)} - X_{\pi(2)}} \\
          & \indent
          - \frac 1{\exp\left(X_{\pi(3)} - X_{\pi(1)}\right) - 1}
          + \frac 1{X_{\pi(3)} - X_{\pi(1)}}
        \end{aligned}
      \right)
      }
    \\
    \shoveleft{
      + \left(
        \frac 1{\exp\left(X_{\pi(3)} - X_{\pi(1)}\right) - 1}
        - \frac 1{X_{\pi(3)} - X_{\pi(1)}} + \frac 1 2
      \right)
      }
    \\
    \shoveleft{
      \cdot
      \left(
        \frac 1{\exp\left(X_{\pi(3)} - X_{\pi(2)}\right) - 1}
        - \frac 1{X_{\pi(3)} - X_{\pi(2)}} + \frac 1 2
      \right)
      \Biggr) = \frac 1 2.
      }
    \\
  \end{multline*}
  This proves the lemma because there is a well-defined $\Q$-linear map
  \begin{gather*}
    \Sym^3_\Q((x_n)_{n \in \N_0}) \to \Q[X_1, X_2, X_3],
    x_i x_j x_k \mapsto \sum_{\pi \in \mathfrak S_3}
    X_{\pi(i)}^i X_{\pi(2)}^j X_{\pi(3)}^k,
  \end{gather*}
  which is injective.
\end{proof}

\begin{proof}[Proof of the theorem.]
  Since $\partial$ is a linear operator of degree $-2$ on $\hBB'$ and
  $\Omega(\mu) = \sum_{k = 0}^\infty \Omega_{2k} \mu^{2k}$, we can
  assume that $\mu = 1$. Set $\Gamma := \sum_{k = 1}^\infty b_{2k}
  w_{2k}$. We have
  \begin{align*}
    \partial\Omega & = \partial \exp(\Gamma)
    = \sum_{n = 0}^\infty \frac{\partial(\Gamma^n)}{n!} \\
    & = \sum_{n = 1}^\infty \partial(\Gamma) \frac{\Gamma^{n - 1}}{(n - 1)!}
    + \sum_{n = 2}^\infty \frac{\partial(\Gamma, \Gamma)} 2 \frac{\Gamma^{n
        - 2}}{(n - 2)!} \\
    & = \left(\partial\Gamma + \frac 1 2 \partial(\Gamma,
      \Gamma)\right) \Omega,
  \end{align*}
  so we see that $\Omega$ is an eigenvector. We just have to
  calculate the eigenvalue:
  \begin{align*}
    \partial\Gamma + \frac 1 2 \partial(\Gamma, \Gamma)
    & = \sum_{k = 1}^\infty b_{2k} \partial w_{2k}
    + \sum_{i, j = 1}^\infty b_{2i} b_{2j} \frac{\partial(w_{2i},
    w_{2j})} 2 \\
    & = \sum_{k = 1}^\infty k b_{2k} \sum_{n = 0}^{2k - 2} w_{n, 2k -
      2 - n}
    + \sum_{i, j = 1}^\infty 4 i j b_{2i} b_{2j} w_{2i - 1, 2j - 1} \\
    & = \frac 1 4 \sum_{k = 1}^\infty \frac{B_{2k}}{(2k)!} \sum_{n =
      0}^{2k - 2} w_{n, 2k - 2 - n}
    + \frac 1 4 \sum_{i, j = 1}^\infty \frac{B_{2i}}{(2i)!}
    \frac{B_{2j}}{(2j)!} w_{2i - 1, 2j - 1} \\
    & = \frac 1 4 \sum_{k = 2}^\infty \frac{B_k}{k!} \sum_{n =
      0}^{k - 2} w_{n, k - 2 - n}
    + \frac 1 4 \sum_{i, j = 2}^\infty \frac{B_i}{i!}
    \frac{B_j}{j!} w_{i - 1, j - 1}.
  \end{align*}

  Applying the (injectiv) map $P$ of lemma~\ref{lem:chmutov} and
  using lemma~\ref{lem:bernoulli} yields
  \begin{multline*}
    P\left(\partial\Gamma + \frac 1 2 \partial(\Gamma, \Gamma)\right)
    = \frac 1 4 P \left(
      \sum_{k = 2}^\infty \frac{B_k}{k!} \sum_{n = 0}^{k - 2}
      w_{n, k - 2 - n}
      + \sum_{i, j = 2}^\infty \frac{B_i}{i!} \frac{B_j}{j!} w_{i - 1,
        j - 1}
    \right) \\
    = \frac 1 {24} x_0^3 = P\left(\frac 1 {48} w_{0, 0}\right)
    =  P\left(\frac 1 {48} \ThetaGraph\right).
  \end{multline*}
  Because of the injectivity of $P$, this proves the theorem.
\end{proof}


\section{Rozansky-Witten classes}

The idea to associate to every graph $\Gamma$ and every hyperk\"ahler
manifold $X$ a cohomology class $\RW_X(\Gamma)$ is due to Rozansky and
Witten (c.f.~\cite{rozansky97}). Kapranov
showed in~\cite{kapranov99} that the metric structure of a
hyperk\"ahler manifold is not nessessary to define these classes. It
was his idea to build the whole theory upon the Atiyah class and the
symplectic structure of a symplectic K\"ahler manifold. We will use
his definition of Rozansky-Witten classes in this section.

\subsection{Rozansky-Witten classes in general}
Let $k$ be a field of characteristic zero, $V$ a finite-dimensional
$k$-vector space, $A = \bigoplus_{i = 0}^\infty A_i$ a $\Z$-graded
(super\-)commutative $k$-algebra and $\sigma$ a symplectic form on $V$.

We shall use this general setting later in the case that $V = \Tang_{X,
  x}$ is the holomorphic tangent space of a complex manifold $X$ at a
point $x$ and $A = \overline{\Omega}^*_{X, x}$ is the Gra\ss mann
algebra of anti-holomorphic forms at $x$.

For every Jacobi diagram $\Gamma$ with $k$ trivalent and $l$ univalent
vertices and every $\alpha \in \Sym^3 V \otimes_k A_1$, we define an
element
\begin{gather}
  \RW_{\sigma, \alpha}(\Gamma) \in \bigwedge^l\nolimits (V^*) \otimes A_k
\end{gather}
by the following procedure:

Let $T$ denote the set of trivalent vertices, $U$ the set of univalent
vertices, $E$ the set of edges and $F$ the set of flags of $\Gamma$.
So, $l = \card U$ and $k = \card T$. The Jacobi diagram $\Gamma$
defines the map
\begin{gather}
  \label{equ:phigamma}
  \begin{aligned}
    \Phi^\Gamma: (\Sym^3 V \otimes A_1)^{\otimes T} \otimes (\End V)^{\otimes
      U}
    & \overset{\text{(1)}}{\To} \left(\bigotimes_{t \in T} V^{\otimes
        t}\right) \otimes A_1^{\otimes T}
    \otimes V^{\otimes U} \otimes (V^*)^{\otimes U} \\
    & \overset{\text{(2)}}{\To}
    V^{\otimes F } \otimes (V^*)^{\otimes U} \otimes A_1^{\otimes T} \\
    & \overset{\text{(3)}}{\To}
    \bigotimes_{e \in E} V^{\otimes e} \otimes
    (V^*)^{\otimes U} \otimes
    A_1^{\otimes T} \\
    & \overset{\text{(4)}}{\To}
    (V^{\otimes 2})^{\otimes\frac{3k + l} 2} \otimes (V^*)^{\otimes l}
    \otimes (A_1)^{\otimes k} \\
    & \overset{\text{(5)}}{\To}
    \bigwedge^l\nolimits (V^*) \otimes A_k,
  \end{aligned}
\end{gather}
where
\begin{enumerate}
\item is induced by the inclusions of the symmetric tensors $\Sym^3
  V$ in the spaces $V^{\otimes t}$ with $t \in T$ (note again that $t$
  is a set of three elements) and the canonical identification $\End V
  = V \otimes V^*$,
\item is induced by the canonical bijection $U \amalg \coprod_{t \in
    T} t \to F$ which is on each $t \in T$ the inclusion of the subset
  $t$ in $F$ and which maps each $u \in U$ to the flag belonging to
  $u$,
\item is induced by the associativity of the tensor product (note that
  $\coprod_{e \in E} e = F$),
\item is induced by choosing total orderings of the set $E$, the sets $e
  \in E$, the set $U$ and the set $T$, which are compatible with the
  orientation of the graph $\Gamma$, and
\item is given by
  $((v_1 \otimes v_2) \otimes \dots \otimes (v_{3k + l - 1}
  \otimes v_{3k + l})) \otimes (\alpha_1 \otimes \dots \otimes
  \alpha_l) \otimes (a_1 \otimes \dots \otimes a_k) \\
  \mapsto
  \sigma(v_1, v_2) \cdot \dots \cdot \sigma(v_{3k + l - 1}, v_{3k + l})
  \cdot (\alpha_1 \wedge \dots \wedge \alpha_l)
  \otimes (a_1 \cdot \dots \cdot a_k)$.
\end{enumerate}
The map $\Phi^\Gamma$ is independent of the specific choice made in (4).
One defines
\begin{equation}
  \RW_{\sigma, \alpha}(\Gamma) := \Phi^\Gamma(\alpha^{\otimes T} \otimes
  (\id_V)^{\otimes U}).
\end{equation}

The following proposition summarizes some of the properties of the
maps $\Phi^\Gamma$ we need later on. All of them follow directly from
the definitions.
\begin{prop}
  Identifying the set of the univalent vertices of $\ell$ with
  the set $\set{1, 2}$, we have
  \begin{gather}
    \label{equ:phiell}
    \Phi^\ell: (\End V)^{\otimes 2} \to \bigwedge^2\nolimits (V^*),
    (v_1 \otimes \lambda_1) \otimes (v_2 \otimes \lambda_2) \mapsto
    \sigma(v_1, v_2) \cdot (\lambda_1 \wedge \lambda_2)
  \end{gather}
  for all $v_1, v_2 \in V$, $\lambda_1, \lambda_2 \in V^*$.
  
  Let $k \in \N_0$. Identifying both the set of the univalent
  vertices and the set of the trivalent vertices of the wheel $w_{2k}$
  (viewed as a Jacobi diagram) with the set $\set{1, \dots, 2k}$, we
  have
  \begin{multline}
    \label{equ:phiwheel}
    \Phi^{w_{2k}}: (\Sym^3 V \otimes A_1)^{\otimes 2k}
    \otimes (\End V)^{\otimes 2k} \to \bigwedge^{2k}\nolimits(V^*)
    \otimes A_{2k}, \\
    \bigotimes_{i = 1}^{2k} (v_i^3 \otimes a_i) \otimes
    \bigotimes_{i = 1}^{2k} (w_i \otimes \lambda_i)
    \mapsto
    - \prod_{i \in \Z/(2k)} (\sigma(v_i, v_{i + 1}) \cdot
    \sigma(v_i, w_i)) \cdot \bigwedge_{i = 1}^{2k} \lambda_i \otimes
    \prod_{i = 1}^{2k} a_i
  \end{multline}
  for all $v_i, w_i \in V$, $\lambda_i \in V^*$ and $a_i \in A_1$.

  Let $\Gamma$ and $\Gamma'$ be two Jacobi diagrams with univalent
  vertices $U$ and $U'$ and trivalent vertices $T$ and $T'$. Let the
  cardinalities of $T$, $T'$, $U$, $U'$ be $k$, $k'$, $l$, $l'$.
  Then the diagram
  \begin{gather}
    \label{equ:phicup}
    \begin{CD}
      \begin{aligned}
        & \left((\Sym^3 V \otimes A_1)^{\otimes T} \otimes (\End
          V)^{\otimes U}\right) \\
        & \indent \otimes
        \left((\Sym^3 V \otimes A_1)^{\otimes T'} \otimes (\End
          V)^{\otimes U'}\right)
      \end{aligned}
      @>{\Phi^\Gamma \otimes \Phi^{\Gamma'}}>>
      \begin{aligned}
        & \left(\bigwedge^{2l}\nolimits(V^*) \otimes A_{2k}\right) \\
        & \indent \otimes
        \left(\bigwedge^{2l'}\nolimits(V^*) \otimes A_{2k'}\right)
      \end{aligned} \\
      @VVV @VVV \\
      (\Sym^3 V \otimes A_1)^{\otimes (T \amalg T')} \otimes (\End
      V)^{\otimes (U \amalg U')}
      @>>{\Phi^{\Gamma \cup \Gamma'}}>
      \bigwedge^{2(l + l')}\nolimits(V^*) \otimes A_{2(k + k')},
    \end{CD} \\
    \intertext{where the vertical maps are the canonical ones, and the
    diagram}
    \label{equ:phiscp}
    \begin{CD}
      \begin{aligned}
        & \left((\Sym^3 V \otimes A_1)^{\otimes T} \otimes (\End
          V)^{\otimes U}\right) \\
        & \indent \otimes
        \left((\Sym^3 V \otimes A_1)^{\otimes T'} \otimes (\End
          V)^{\otimes U'}\right)
      \end{aligned}
      @>{\Phi^\Gamma \otimes \Phi^{\Gamma'}}>>
      \begin{aligned}
        & \left(\bigwedge^{2l}\nolimits(V^*) \otimes A_{2k}\right) \\
        & \indent \otimes
        \left(\bigwedge^{2l'}\nolimits(V^*) \otimes A_{2k'}\right)
      \end{aligned} \\
      @AAA @VVV \\
      (\Sym^3 V \otimes A_1)^{\otimes T} \otimes
      (\Sym^3 V \otimes A_1)^{\otimes T'}
      @>>{\Phi^{\left<\Gamma, \Gamma'\right>}}>
      A_{2(k + k')},
    \end{CD}
  \end{gather}
  where the left vertical map is induced by tensoring with the
  identity on $V$ and the right one is just~\eqref{equ:formscp},
  commute.
\end{prop}

\subsection{Rozansky-Witten invariants of symplectic K\"ahler
  manifolds}

All our K\"ahler manifolds are assumed to be compact.

Let $X$ be a symplectic K\"ahler manifold, $\sigma$ a fixed symplectic
form and $\alpha$ a Dolbeault representative of the Atiyah class of
$X$. We can use $\sigma$ to identify the tangent bundle $\TangX$ of
$X$ with the cotangent bundle $\Omega_X$. Doing this, $\alpha$ can be
viewed as an element of $A^1(X, \TangX^{\otimes 3})$.  Here, $A^k(X,
E)$ denotes the space of $(0, k)$-forms with values in a holomorphic
vector bundle $E$. We set $A^{l, k}(X) := A^k(X,
\Omega_X^l)$.

The following proposition was proven by
Kapranov in~\cite{kapranov99}.
\begin{prop}
  The values of $\alpha$ are symmetric tensors, i.e.\ $\alpha \in
  A^1(X, \Sym^3\Tang_X)$.
\end{prop}

\begin{defn}
  For every Jacobi diagram $\Gamma$ with $k$ trivalent and $l$
  univalent vertices, one defines
  \begin{gather}
    \RW_{\sigma}(\Gamma) \in \HH^k(X, \Omega^l)
  \end{gather}
  to be the Dolbeault cohomology class of the ($\bar\partial$-)closed
  $(l, k)$-form
  \begin{gather}
    \label{equ:rwdef}
    (x \mapsto \RW_{\sigma_x, \alpha_x}(\Gamma)) \in A^{l, k}(X).
  \end{gather}
  (That the form defined in~\eqref{equ:rwdef} is
  $\bar\partial$-closed follows from the fact that $\sigma$ and
  $\alpha$ are $\bar\partial$-closed.)

  For a $\Q$-linear combination $\gamma$ of Jacobi diagrams,
  $\RW_\sigma(\gamma)$ is defined by linear extension.
\end{defn}

In~\cite{kapranov99}, Kapranov also showed the following proposition, which
is crucial for the next definition.
\begin{prop}
  If $\gamma$ is a $\Q$-linear combination of Jacobi diagrams that is
  zero modulo the anti-symmetry and IHX relations, then
  $\RW_\sigma(\gamma) = 0$.
\end{prop}

\begin{defn}
  We define a linear map
  \begin{gather}
    \RW_\sigma: \hBB \to \HH^*(X, \Omega^*),
  \end{gather}
  which maps $\hBB_{k, l}$ into $\HH^k(X, \Omega^l)$, by mapping a
  homology class of a Jacobi diagram $\Gamma$ to $\RW_{\sigma}(\Gamma)$.

  The values of the just defined map $\RW_\sigma$ are called
  \emph{Rozansky-Witten classes of the symplectic K\"ahler manifold $X$}.
\end{defn}

\subsection{Examples of Rozansky-Witten classes}

Let $X$ and $\sigma$ be as before.

\begin{exmp}
  The Dolbeault cohomology class $[\sigma] \in \HH^{2, 0}(X)$ is a
  Rozansky-Witten class, more precisely, we have
  \begin{gather}
    \RW_\sigma(\ell) = 2 \sigma,
  \end{gather}
  which follows from~\eqref{equ:phiell}.
\end{exmp}

The following example is due to Hitchin and Sawon~\cite{hitchin99}. It
is of great importance for their and our results. 
\begin{exmp}
  Let $\ch(X) = \sum_{k = 0}^\infty s_{2k}/(2k)!$, $s_{2k} \in
  \HH^{2k, 2k}(X)$, be the Chern character of $X$. Then
  \begin{gather}
    \label{equ:rwwheel}
    \RW_\sigma(w_{2k}) = - s_{2k}
  \end{gather}
  for all $k \in \N_0$. (Note that for a symplectic manifold $\ch(X)$
  has no term in degree $(k, k)$ for $k$ odd.)

  Since the algebra of characteristic classes of $X$ is spanned by the
  classes $s_{2k}$, every characteristic class is a Rozansky-Witten
  class due to proposition~\ref{prop:rwalg} below.
\end{exmp}

A proof of \eqref{equ:rwwheel} is given by Hitchin and Sawon
in~\cite{hitchin99}, where the Rozansky-Witten invariants are
defined by using the Riemann curvature tensor of a hyperk\"ahler
metric of $X$ instead of the Atiyah class. An idea of the proof in our
context is given below.  

\begin{proof}[Idea of the proof of \eqref{equ:rwwheel}]
  As remarked in~\cite{kapranov99} by Kapranov, a Dolbeault
  representative of the Atiyah class $\alpha \in A^1(X, \Omega_X
  \otimes \End\TangX)$ on every K\"ahler manifold $X$ is related to
  the characteristic classes $s_k \in \HH^{k, k}(X)$ of $X$ in the
  following way: A Dolbeault representative of $s_k$ is given by
  $\Alt(\tr(\alpha^k))$. Here, $\alpha^k \in A^k(X, \Omega_X^{\otimes
    k} \otimes \End\TangX)$ means taking the $k$-th product of
  $\alpha$ viewed as $(0, 1)$-form, giving an element of $A^k(X,
  \Omega_X^{\otimes k} \otimes_X (\End\TangX)^{\otimes k})$, and then
  using the associative algebra structure of $\End\TangX$.
  Further, $\tr$ means taking the trace on $\End\TangX$, and
  $\Alt$ is induced by the canonical projection $\Omega^{\otimes k}_X
  \to \Omega^k_X$.
  
  That this procedure on symplectic K\"ahler manifolds is essentially
  the same as taking the Rozansky-Witten class of a wheel follows from
  ~\eqref{equ:phiwheel}.
\end{proof}

\begin{exmp}
  The Todd genus of a symplectic complex manifold $X$ is given by
  \begin{gather}
    \td(X) = \exp\left(- 2 \sum_{k = 0}^\infty b_{2k} s_{2k}\right),
  \end{gather}
  with $b_{2k}$ being a modified Bernoulli number
  (see~\cite{hitchin99} for this statement). Thus,
  \begin{gather}
    \td(X) = \RW_\sigma(\Omega^2).
  \end{gather}
\end{exmp}

\subsection{Properties of Rozansky-Witten classes}

\begin{prop}
  \label{prop:rwalg}
  The map $\RW_\sigma: \hBB \to \HH^{*, *}(X)$ is a morphism of
  graded $\Q$-algebras.

  Furthermore, we have
  \begin{gather}
    \RW_\sigma(\left<\Gamma, \Gamma'\right>)
    = \left<\RW_\sigma(\Gamma), \RW_\sigma(\Gamma')\right>
  \end{gather}
  for all $\Gamma, \Gamma' \in \hBB$.
\end{prop}

\begin{proof}
  The statements follow from~\eqref{equ:phicup} and~\eqref{equ:phiscp}.
\end{proof}

The following proposition is also stated in~\cite{hitchin99} in a slightly
different notation. With the formalism introduced here, we can give a
compact proof.

\begin{prop}
  If $X$ is irreducible, we have
  \begin{gather}
    \RW_\sigma(\ThetaGraph) = \frac{2 \int_X c_2(X) \exp(\sigma + \bar\sigma)}
    {n \int \exp_X(\sigma + \bar\sigma)} \cdot [\bar\sigma].
  \end{gather}
\end{prop}

\begin{proof}
  Due to the irreducibility of $X$, i.e.\ $\HH^{2k}(X, \OX) = \C \cdot
  [\bar\sigma]$ for all $k \in \N_0$, we can write
  \begin{gather}
    \alpha = \frac{\int_X \alpha \exp(\sigma + \bar\sigma)}{n \int_X \exp
    (\sigma + \bar\sigma)} \cdot [\bar\sigma]
  \end{gather}
  for all $\alpha \in \HH^{0, 2}(X)$.
  Using this, we have
  \begin{multline}
    \RW_\sigma(\ThetaGraph) = \RW_\sigma\left(\frac 1 2 \left<w_2,
      \ell\right>\right)
    = \frac 1 2 \left<\RW_\sigma(w_2), \RW_\sigma(\ell)\right> \\
    = \left<-s_2, \sigma\right> = 2 \left<c_2(X), \exp \sigma\right>
    = \frac {2 \int_X \left<c_2(X), \exp\sigma\right> \exp(\sigma +
      \bar\sigma)}
    {n \int_X \exp(\sigma + \bar\sigma)} \cdot [\bar\sigma],
  \end{multline}
  which proves the proposition because of proposition~\ref{prop:laexp}.
\end{proof}


\section{The Euler characteristic of a line bundle in terms of the
quadratic form}

Let $X$ be an irreducible symplectic K\"ahler manifold with symplectic
form $\sigma$. 

\subsection{The quadratic form on an irreducible symplectic manifold}

Let us assume $\int_X (\sigma\bar\sigma)^n = 1$ for this subsection.

\begin{defn}
  We set
  \begin{multline}
    f_X: \HH^2(X, \C) \to \C, \\
    \alpha \mapsto
    \frac n 2 \int_X \alpha^2 (\sigma\bar\sigma)^{n - 1}
    + (1 - n)
    \left(\int_X \alpha \bar\sigma (\sigma\bar\sigma)^{n - 1}\right)
    \left(\int_X \alpha \sigma (\sigma\bar\sigma)^{n - 1}\right).
  \end{multline}
\end{defn}

\begin{rmk}
  The form $f_X$ is a multiple of the Beauville-Bogomolov quadratic form $q_X$
  on $X$. A detailed discussion of $f_X$ and $q_X$ can be found
  in~\cite{huybrechts99}. Note that $f_X(\alpha) = \frac n 2 \int
  \alpha^2(\sigma\bar\sigma)^{n - 1}$ if $\alpha \in \HH^{1, 1}(X)$.
\end{rmk}

We need a fact about $f_X$, which was proved by Huybrechts
in~\cite{huybrechts99}.

\begin{prop}
  \label{prop:daniel}
  Assume that $\alpha \in \HH^{4j}(X, \C)$ is a characteristic class
  of $X$ (it suffices to assume that $\alpha$ is of type $(2j, 2j)$ on
  all small deformations of $X$). Then there is a constant $c_{\alpha}
  \in \C$, depending on $\alpha$, with
  \begin{gather}
    \label{equ:daniel}
    \int_X \alpha \beta^{2(n - j)} = c_\alpha \cdot (f_X(\beta))^{n - j}
  \end{gather}
  for all $\beta \in \HH^2(X, \C)$.
\end{prop}

\begin{cor}
  \label{cor:daniel}
  For all $\beta \in \HH^2(X, \C)$, we have
  \begin{equation}
    \int_X \beta^{2n} = \binom{2n} n (f_X(\beta))^n
  \end{equation}
  and
  \begin{equation}
    \int_X c_2(X) \beta^{2n - 2} = \binom{2n - 2} {n - 1} 
    \left(\int_X c_2(X) (\sigma\bar\sigma)^{n - 1}\right) \cdot
    (f_X(\beta))^{n - 1}.
  \end{equation}
\end{cor}

\begin{proof}
  Set $\alpha = 1$ in~\eqref{equ:daniel}. To obtain $c_1$, compute both
  sides with $\beta = \sigma + \bar\sigma$, which proves the first equation.
  
  Then set $\alpha = c_2(X)$ in~\eqref{equ:daniel}. To obtain $c_{c_2(X)}$
  compute both sides with $\beta = \sigma + \bar\sigma$, which leads
  to the second equation.
\end{proof}

\subsection{A Hirzebruch-Riemann-Roch formula}
\label{s:hrrformula}
Let us define the number
\begin{gather}
  \lambda := \frac{24 n \int_X \exp(\sigma + \bar\sigma)}{\int_X
    c_2(X) \exp(\sigma + \bar\sigma)}
\end{gather}
for the pair $(X, \sigma)$.  It is well-defined (i.e.\ the denominator
does not vanish) because $\int_X c_2(X) \exp(\sigma + \bar\sigma)$
equals the $\mathrm L^2$-norm of the Riemann curvature tensor of $X$
(having been equipped with a hyperk\"ahler metric compatible with the
given symplectic structure) up to a positive constant
(see~\cite{hitchin99}). But if the Riemann curvature tensor vanishes,
$X$ is a torus, which contradicts the assumption on irreducibility.

We continue to write $b_{2k}$ for the modified Bernoulli
numbers and $s_{2k}$ for the homogeneous components of $\ch(X)$
(up to the constant $(2k)!$).

\begin{prop}
  \label{prop:todd}
  For all $\mu \in \C$ with $\mu^2 + \mu^{-2} = 2 + \lambda$, we have
  \begin{gather}
    \label{equ:todd}
    \int_X \td(X) \exp(\sigma + \bar\sigma) =
    \int_X \exp\left(- \sum_{k = 0}^\infty
      b_{2k}s_{2k}\cdot(\mu^{2k} + \mu^{-2k})\right).
  \end{gather}
\end{prop}

\begin{proof}
  Note that $\left<\Omega, \Omega\right> = \left(\Omega(\mu),
  \Omega(\mu^{-1})\right>$. Then
  \begin{multline*}
    \int_X \td(X) \exp(\sigma + \bar\sigma)
    = \int_X \RW_\sigma(\Omega^2) \exp(\sigma + \bar\sigma)
    \\
    \shoveleft{
      = \int_X \left<\RW_\sigma(\Omega^2),
        \exp\sigma\right>\exp(\sigma + \bar\sigma)
      = \int_X \RW_\sigma\left(\left<\Omega^2, \exp(\ell/2)\right>\right)
      \exp(\sigma + \bar\sigma)
      }
    \\
    \shoveleft{
      = \int_X \RW_\sigma\left(\left<\exp(\partial)(\Omega^2),
          1\right>\right) \exp(\sigma + \bar\sigma) 
      } \\
    \shoveleft{
      = \int_X \RW_\sigma\left(\left<\exp(\partial)\Omega,
          \exp(\partial)\Omega\right>\right) \exp(\sigma +
      \bar\sigma)
      } \\
    \shoveleft{
      =\int_X \RW_\sigma\left(\left<\exp(\partial)\Omega,
          \exp(\partial)\Omega\right>\right)
      \RW_\sigma\left(\exp\left(\frac\lambda{48}
          \ThetaGraph\right)\right) \exp(\sigma) 
      } \\
    \shoveleft{
      = \int_X \RW_\sigma\left(\exp\left(\frac{2 + \lambda}{48}
          \ThetaGraph \right)
        \left<\Omega, \Omega\right>\right) \exp(\sigma) 
      } \\
    \shoveleft{
      = \int_X \RW_\sigma\left(\exp\left(\frac{2 + \lambda}{48}
          \ThetaGraph \right)
        \left<\Omega(\mu), \Omega(\mu^{-1})\right>\right)
      \exp(\sigma)
      } \\
    \shoveleft{
      = \int_X \RW_\sigma\left(\left<\exp(\partial)\Omega(\mu),
          \exp(\partial)\Omega(\mu^{-1})\right>\right) \exp(\sigma)
      } \\
    \shoveleft{
      = \int_X
      \RW_\sigma\left(\left<\exp(\partial)(\Omega(\mu)\Omega(\mu^{-1})),
          1\right>\right) \exp(\sigma)
      } \\
    \shoveleft{
      = \int_X \RW_\sigma(\left<\Omega(\mu)\Omega(\mu^{-1}),
        \exp(\ell/2)\right>) \exp(\sigma)
      } \\
    \shoveleft{
      = \int_X \left<\RW_\sigma(\Omega(\mu)\Omega(\mu^{-1})),
        \exp(\sigma)\right>\exp(\sigma)
      = \int_X \RW_\sigma(\Omega(\mu)\Omega(\mu^{-1})) \exp(\sigma)
      }
    \\
    \shoveleft{
      = \int_X \RW_\sigma(\Omega(\mu)\Omega(\mu^{-1}))
      = \int_X \exp\left(-\sum_{k =
          0}^\infty b_{2k}s_{2k} \cdot (\mu^{2k} + \mu^{-2k})\right).
      }
    \\
  \end{multline*}
\end{proof}

\begin{rmk}
  \label{rmk:todd}
  The right hand side of~\eqref{equ:todd} is polynomial in $\lambda$
  since
  \begin{equation}
    \label{equ:chebychev}
    \mu^n + \mu^{-n} = 2 T_n \left(\pm \sqrt{\frac z 4 + \frac 1 2}\right)
  \end{equation}
  for all $\mu, z \in \C$ with $\mu^2 + \mu^{-2} = z$ and $n \in
  \N_0$. Here, $T_n$ denotes the $n^{\textrm{th}}$ Chebychev
  polynomial. Note that $T_n$ is even for $n$ even.
\end{rmk}

\begin{proof}
  We can assume that $z \in [-2, 2]$. If follows that $\mu = \pm e^{\pm
  ix}$ with $x = \frac 1 2 \arccos(\frac z 2) =
  \arccos\left(\sqrt{\frac z 4 + \frac 1 2}\right)$. By definition of
  the Chebychev polynomials, we have
  \begin{gather*}
    \mu^n + \mu^{-n} = 2 \cos\left(n \cdot \arccos\left(\pm
        \sqrt{\frac z 4 + \frac 1 2}\right)\right) = 2 T_n\left(\pm
      \sqrt{\frac z 4 + \frac 1 2}\right).
  \end{gather*}
\end{proof}

\begin{defn}
  For every line bundle $L$ on $X$, we call the rational number
  \begin{equation}
    \lambda(L) := \begin{cases}
      \frac{24 n \int_X \ch(L)}{\int c_2(X) \ch(L)} & \text{if
        well-defined} \\
      0 & \text{otherwise}
    \end{cases}
  \end{equation}
  the \emph{characteristic value of $L$}.
\end{defn}

We could have defined a characteristic value $\lambda(\alpha)$ for
every cohomology class $\alpha \in \HH^2(X, \C)$ so that $\lambda(L) =
\lambda(c_1(L))$. Then the number $\lambda$ for the pair $(X, \sigma)$
defined above is just $\lambda(\sigma + \bar\sigma)$.

\begin{rmk}
  If $\lambda(L) = 0$, then
  \begin{equation}
    \chi(L) = \chi(\OX) = \int_X \td(X).
  \end{equation}
\end{rmk}

\begin{proof}
  If $\lambda(X) = 0$, then $\int_X \ch(L) = 0$ or $\int_X c_2(X)
  \ch(L) = 0$. In both cases, $f_X(c_1(L)) = 0$ due to the
  corollaries~\ref{cor:daniel} and~\ref{cor:daniel}. Because of the
  proposition~\ref{prop:daniel} and the usual Hirzebruch-Riemann-Roch
  formula, the assertion follows.
\end{proof}

\begin{prop}
  \label{prop:lambda}
  For the characteristic value of line bundle $L$ on $X$,
  \begin{equation}
    \lambda(L) = \frac{12 \int_X c_1(L)^2 (\sigma\bar\sigma)^{n-1}}{\int
      c_2(X) (\sigma\bar\sigma)^{n - 1}}.
  \end{equation}
  Thus, $\lambda(L)$ is up to a multiple the value of the quadratic
  form $f_X$ applied to the first Chern class of $L$.
\end{prop}

\begin{proof}
  This is an immediate application of corollary~\ref{cor:daniel}.
\end{proof}

\begin{thm}
  Let $X$ be an irreducible symplectic K\"ahler manifold with
  symplectic form $\sigma$.

  For every line bundle $L$ on $X$, the Euler characteristic of $L$
  can be expressed as
  \begin{equation}
    \label{equ:hrr}
    \chi(L) = \int_X \exp\left(-2 \sum_{k = 1}^\infty b_{2k} s_{2k}
      T_{2k}\left(\sqrt{\lambda(L)/4 + 1}\right)\right).
  \end{equation}

  The right hand side is polynomial in $\lambda(L)$ and thus
  polynomial in the quadratic form of $c_1(L)$.
\end{thm}

\begin{proof}
  We can assume that $\lambda(L) \neq 0$.
  Since both sides of \eqref{equ:todd} are rational functions in the cohomology
  class of $\sigma + \bar\sigma$ and the equality holds for every irreducible
  symplectic complex manifold, by a standard argument of
  hyperk\"ahler geometry, \eqref{equ:todd} remains true, if we
  substitute $\sigma + \bar\sigma$ by an arbitrary cohomology class
  in $\HH^2(X)$. Here, we have to substitute $\sigma + \bar\sigma$ by
  $c_1(L)$.
\end{proof}

\begin{example}
  Let $\lambda := \lambda(L)$ be the characteristic value of $L$. Then
  \begin{equation}
    \begin{split}
      & \td(\TangX(L)) = 1
      + \frac 1 {12} \left(c_2 + \frac 1 2 c_2 \cdot \lambda\right) \\
      & + \frac 1 {720} \left(3 c_2^2 - c_4 
        + \left(\frac 7 2 c_2^2 - 2 c_4\right)
        \lambda + \left(\frac 7 8 c_2^2 - \frac 1 2 c_4\right)
        \lambda^2\right) \\
      & + \frac 1 {30240} \left(
        \begin{aligned}
          & \left(5 c_2^3 - \frac 9 2 c_2 c_4 + c_6\right)
          + \left(\frac{41} 4 c_2^3 - \frac{53}4 c_2 c_4 + \frac 9 2
            c_6\right) \lambda \\
          & \indent + \left(\frac{93}{16} c_2^3 - \frac{33} 4 c_2 c_4 + 3
            c_6\right) \lambda^2
          + \left(\frac{31}{32} c_2^3 - \frac{11} 8 c_2 c_4 + \frac 1 2
            c_6\right) \lambda^3
        \end{aligned}
      \right) \\
      & + \cdots,
    \end{split}
  \end{equation}
  where the $c_i$ denote the Chern classes of $X$.
\end{example}

There are (at least) two possibilities to write~\eqref{equ:hrr} in
another form. One is to define ``deformed Todd classes''
$\td_\epsilon$, $\epsilon \in \C$, for arbitrary complex manifolds $X$:

\begin{defn}
  Let $X$ be an arbitrary complex manifold. We set
  \begin{equation}
    \td_\epsilon(X) = \exp\left(\sum_{k = 0}^\infty t_k
      T_k(\epsilon)\right),
  \end{equation}
  where $t_k \in \HH^{k, k}$ such that $\ln\td(X) = \sum_{k =
    0}^\infty t_k$, i.e.\ $\td_1(X) = \td(X)$.
\end{defn}

Applying this definition to our irreducible symplectic K\"ahler
manifold $X$, the equation~\eqref{equ:hrr} becomes
\begin{equation}
  \chi(L) = \int_X \td_\epsilon(X)
\end{equation}
with $\epsilon = \sqrt{\lambda(L)/4 + 1}$.

Let us recall some facts about the
Grothendieck ring $\KK^0(X)$ of complex vector bundles over $X$ (see
for example~\cite{atiyah67}). After tensoring with $\C$, the Chern
character gives us a ring homomorphism $\ch: \KK^0(X, \C) \to \HH^*(X,
\C)$ and taking the Todd class gives us a group homomorphism $\td:
\KK^0(X, \C) \to 1 + \HH^{* > 0}(X, \C)$. Furthermore, there are the Adams
operations $\psi^p: \KK^0(X, \C) \to \KK^0(X, \C)$, $p \in \N$, which
are ring homomorphisms with $\psi^p(L) = L^p$ for
every line bundle $L$. The $\psi^p$ commute, and we can write
$\KK^0(X, \C) = \bigoplus_{k = 0}^\infty \Gr^{2k} \KK^0(X, \C)$ such that
$\Gr^{2k} \KK^0(X, \C)$ is the eigenspace of $\psi^p$ for $p \geq 2$ to
the eigenvalue $p^k$. Using this grading on $\KK^0(X, \C)$, the Chern
character becomes a homomorphism of graded rings.

If $M = \bigoplus_{i \in \N_0} M_i$ is a graded $R$-module, $R$ a
commutative ring and $(\lambda_i)_{i \in \N_0}$ is a sequence in $R$,
we define a morphism of graded modules $(\lambda_i)_{i \in \N_0} \cdot:
M \to M, m \in M_i \mapsto \lambda_i \cdot m$. This definition will be
applied to $\KK^0(X, \C)$:

\begin{defn}
  For every complex number $\epsilon \in \C$ and complex manifold $X$, one
  defines a homomorphism
  \begin{equation}
    \phi_\epsilon: \KK^0(X, \C) \to \KK^0(X, \C), E \mapsto
    \left(T_i(\epsilon)\right)_{i \in \N_0} \cdot E
  \end{equation}
\end{defn}

Applying this definition to our irreducible symplectic K\"ahler
manifold $X$, the equation~\eqref{equ:hrr} becomes
\begin{equation}
  \chi(L) = \int_X \td(\phi_\epsilon(\TangX))
\end{equation}
with $\epsilon = \sqrt{\lambda(L)/4 + 1}$.


\section{Generalized Kummer Varieties.}

\subsection{Basic Facts}

Generalized Kummer varieties were introduced by Beauville
in~\cite{beauville83} as the second series of examples of
irreducible symplectic manifolds. We briefly recall their
construction and some basic properties of these varieties.

We start with an abelian surface $A$. Let $\SA{n}:= A^n/\mf{S}_n$
be the $n$-fold symmetric product, denote by $\HA{n}$ the $n$-th
Hilbert scheme of points on $A$, and let $\rho:\SA{n}\To \HA{n}$
be the Hilbert-Chow-morphism. Since summation $A^n\To A$ is
symmetric, it factors over $\SA{n}$, and by composing with $\rho$
one gets the ``summation'' morphism $s:\HA{n}\To A$. Now the
$(n-1)$-th Kummer variety $\KA{n-1}$ is defined as the fiber over
$0$ of the summation morphism $s$. As was shown by Beauville,
$\KA{n-1}$ is an irreducible symplectic variety of dimension
$2(n-1)$. In order to fix notations, we briefly show that the
$(n-1)$-th Kummer variety is smooth.

This can be seen as follows: Denote the translation operation on
$\HA{n}$ by an element $a\in A$ by $t_a$. Considering $A$ acting
on $\HA{n}$ via $t_a$ while acting on $A$ via $t_{na}$, one sees
that $s$ is an equivariant morphism. Since $A$ acts transitively
on itself, all fibers are isomorphic. (In particular, the
definition of $\KA{n-1}$ is independent of the choice of $0\in
A$.) Since there are smooth fibers, the Kummer variety is smooth.

Actually, the fibration $\HA{n}\stackrel{s}{\To} A$ is isotrivial,
i.e. one has the following Cartesian diagram, where $n$ denotes
the morphism ``multiplication by $n$'' $A \xrightarrow{\Delta}
A^n\xrightarrow{\Sigma}A$:

\begin{equation}
  \label{isotriviality}
  \begin{CD}
    K^{n-1}A \times_{\IC} A @>{\nu}>> \HA{n} \\
    @V{p_A}VV @VV{s}V \\
    A @>>n> A.
  \end{CD}
\end{equation}

In terms of closed points, the fibre product is
$\{(\xi,a)\in\HA{n}\times A| s(\xi)=na\}$. This is isomorphic to
$K^{n-1}A\times A$ via $(\xi,a)\mapsto (t_{-a}(\xi),a)$.
Therefore, on closed points the morphism $\nu$ in the above
diagram is given by $(\xi,a)\mapsto t_a(\xi)$.

\begin{example} \emph{The classical Kummer surface}.

The easiest example -- and the reason for the terminology
`generalized Kummer varieties' -- is the classical Kummer surface.
In the classical context this surface is constructed as follows
(for details, cf.\ e.g.~\cite{beauville96}):

One starts with an Abelian surface $A$ and considers the singular
quotient $A/\mspace{-8mu}\sim$ by the involution $(-1)_A$. The
singularities are the images of the 16 two-division points. The
desingularization, which we again denote by $K^1A$, is the
classical Kummer suface which is an K3 surface.

Alternatively one can first blow up the 16 points of order 2 of A
and let $K^1A$ be the quotient of the induced involution on the
blown up surface $\hat{A}$. In other words, one has the following
diagram:

\begin{equation*}
  \begin{CD}
    \KA 1 @<<< \hat A \\
    @V{\epsilon}VV @VVV \\
    A/\mspace{-8mu}\sim @<<< A.
  \end{CD}
\end{equation*}

$\KA 1$ can be identified with the fiber over 0 of $s:\HA{2}\To A$
as follows:

In the case $n=2$ the Hilbert Chow morphism
$\rho_2:\HA{2}\To\SA{2}$ is simply the blow-up of the diagonal
$\Delta\subset \SA{2}$ (cf.~\cite{fogarty68}). Denote by
$\tilde{\Delta}:A\To \SA{2}$ the morphism induced from the two
isomorphisms $\id_A$ and $(-1)_A$. On closed points,
$\tilde{\Delta}$ is given by $a\mapsto (a,-a)$. By definition of
$\SA{2}$ and $A/\mspace{-8mu}\sim$, it descends to an morphism
$A/\mspace{-8mu}\sim \To \SA{2}$, which we again denote by
$\tilde{\Delta}$. From the universal property of the fibre
product, we get an isomorphism $A/\mspace{-8mu}\sim \To
\Sigma^{-1}(0)$.

So, we have the following diagram:

\begin{equation}
  \label{ClassCase}
  \begin{CD}
    \KA 1 @>>> \HA 2 \\
    @V{\epsilon}VV @VV{\rho}V \\
    A/\mspace{-8mu}\sim @>{\tilde\Delta}>> \SA 2 \\
    @VVV @VV{\Sigma}V \\
    0 @>>> A.
  \end{CD}
\end{equation}

This shows that $\eps=\rho|_{\KA 1}$ and that the two descriptions
of $\KA 1$ coincide.
\end{example}

\subsection{Explicit Hirzebruch-Riemann-Roch for $\KA{n}$.}

In this section we will prove the following
Hirzebruch-Riemann-Roch formula formula without using
Rozansky-Witten classes in a purely algebro-geometric way.

\begin{thm}\label{eRR}
Let $L$ be a line bundle on $\KA{n}$. The Euler characteristics of
$L$ is given by
\[
    \chi(L) = (n+1) \binom{\frac{(n+1)}{4}\lambda(L) + n}{n},
\]
where $\lambda(L)$ is the characteristic value of $L$ defined in
section~\ref{s:hrrformula}.
\end{thm}

We will prove the theorem as follows: First, observe that since
the Hirzebruch-Riemann-Roch formula has the form $\chi(L)= \sum
a_2i/(2i)! q_X(L)^i$ with universal coefficients $a_i$, it is
actually enough to consider a single line bundle $L$ with
$\chi(L)\neq 0$ and computing the coefficients by calculating
$\chi(L^m)$. We will prove the theorem using special line bundels
$K^{n}L$ on $\KA{n}$, which are constructed from an invertible
sheaf $L\in \Pic(A)$ as follows:

In a first step define $S^{(n+1)}L := \pi_*{(L^{\boxtimes
{(n+1)}})}^{\mf{S}_{(n+1)}}$ the $\mf{S}_{(n+1)}$-symmetrized line
bundle on $\SA{n+1}$. Here, $\pi:A^{n+1}\to \SA{n+1}$ is the
projection and $L^{\boxtimes n}$ denotes the sheaf
$\bigotimes_{1=1}^n\pr_i^*L $. Now define $L_{(n+1)} :=
\rho_{(n+1)}^* S^{(n+1)}L \in \Pic(\HA{n+1})$. We are interested
in the cohomology of the restricted bundle $K^{n}L :=
L_{(n+1)}|_{K^{n}A}$ on the $n$-th Kummer variety.

One knows that for such bundles the quadratic form
$f_X(c_1(K^{n}L))$ coincides -- up to a positive scalar factor --
with the intersection pairing $c_1(L)^2$ on $\HH^2(A,\IZ)$
(cf.~\cite{beauville83}). If we normalize the quadratic form by
defining $q(K^{n}L):=c_1(L)^2$, our theorem follows from the
following

\begin{lem}\label{RRprop}

The Euler characteristics of $L\in \Pic(A)$ and of $K^{n}L\in
\Pic(\KA{n})$ are related by
\[
    \chi(K^{n}L) = (n+1) \binom{\frac{c_1(L)^2}{2} + n}{n}.
\]
\end{lem}

\begin{proof}[Proof of the Theorem]
Using the lemma above the theorem is proven if we show that
$q=\frac{n+1}{2}\lambda$.

By definition of $K^{n}L$ and $q$ we get
\[
    q((K^{n}L)^{\otimes m}) = q(K^{n}(L^{\otimes m})) =
    m^2q(K^{n}L)
\]
for every $m\in \IN$.

Now HRR in its classical form -- using $c_1(X)=0$ -- gives
\[
    \chi(K^{n}L^{\otimes m}) = m^{2n}\int\ch(K^{n}L) +
    m^{2n-2}\int \frac{c_2(X)}{12}\ch(K^{n}L) + \text{ terms of lower degree,}
\]
whereas our formula in lemma \ref{RRprop} reads
\[
    \chi(K^{n}L^{\otimes m}) = m^{2n}\frac{n+1}{2^n n!}q(K^{n}L)^n + m^{2n-2}
    \frac{(n+1)^2 n}{2^n n!}q(K^{n}L)^{n-1} + {\rm told.}
\]

Comparing coefficients, we find
\[
    \lambda(K^{n}L) =
    \frac{24n\int\ch(K^{n}L)}{\int c_2(X)\ch(K^{n}L)} =
    \frac{2}{n+1}q(K^{n}L).
\]

Since both $\lambda$ and $q$ are positive multiples of $f_X$ and
the coefficients of the HRR formula as a polynomial in the
quadratic form are independent of the line bundle considered, the
theorem follows.
\end{proof}

The rest of this section is dedicated to the proof of lemma
\ref{RRprop}. To simplify notations we will prove the assertion
for $K^{n-1}L \in \Pic(\KA{n-1})$.

In the next lemma we will compute the Euler characteristic of the
symmetrized line bundle $S^nL$.

\begin{lem}
Let $L$ be a line bundle on $A$. Then one has
\[
    \chi(S^nL) = \binom{\chi(L) +n-1}{n}
\]
\end{lem}

\begin{proof}

Let $L, H \in \Pic(A)$. One has
\[
    S^n(L\otimes H) = \pi_*((L\otimes H)^{\boxtimes n})^{\mf{S}_n}
    = \pi_*(L^{\boxtimes n}\otimes H^{\boxtimes n})^{\mf{S}_n},
\]

and since the action of the symmetric group does not flip the
factors $L^{\boxtimes n}$ and $H^{\boxtimes n}$, and for an
arbitrary line bundle $M$ on $A$ we have $\pi^*S^nM = M^{\boxtimes
n}$, we get
\[
    S^n(L\otimes H) = \pi_*(L^{\boxtimes n} \otimes
    \pi^*S^nH)^{\mf{S}_n} = (\pi_*(L^{\boxtimes n})\otimes
    S^nH)^{\mf{S}_n} = S^nL\otimes S^nH
\]

Now, if $H$ is an ample invertible sheaf, $H^{\boxtimes n}$ is
also ample and so is $S^nH$, because $\pi$ is a finite surjective
morphism and $\pi^* S^nH = H^{\boxtimes n}$ is ample. Let $N\in
\IN$ be large enough that both $L\otimes H^N$ and $S^n(L\otimes
H^N) = S^nL\otimes (S^nH)^N$ have no higher cohomology such that
$\chi(S^n(L\otimes H^N)) = h^0(S^n(L\otimes H^N))$.

For the global sections of an symmetrized line bundle $S^nM$, one
has the isomorphism

\begin{equation}\label{H0SnL}
    \HH^0(\SA{n}, S^nM) \simeq \HH^0(A^n, M^{\boxtimes n})^{\mf{S}_n}
    = (\HH^0(A,M)^{\otimes n})^{\mf{S}_n} \simeq \Sym^n \HH^0(A,M).
\end{equation}

Now replacing $M$ by $L\otimes H^N$, we find

\[
    \chi(S^n(L\otimes H^N)) = \binom{\chi(L\otimes H^N) +n-1}{n}\quad \text{ for all } N\gg
    0.
\]

Since $N\mapsto \chi(S^n(L\otimes H^N)) = \chi(S^nL\otimes
(S^nH)^N)$ is a polynomial in $N$, evaluation in $N=0$ proves the
lemma.
\end{proof}

Consider the Hilbert-Chow morphism $\rho:\HA{n}\To \SA{n}$. Since
$\rho$ is a birational proper morphism of normal varieties one has
$\rho_*\shO_{\HA{n}} = \shO_{\SA{n}}$. Furthermore, $\SA{n}$ has
rational singularities as quotient of a smooth variety by a finite
group (cf.~\cite{kollar98}). Therefore its resolution $\rho$
satisfies $R^j\rho_*\shO_{\HA{n}} =0$ for $j>0$. Using the Leray
spectral sequence one gets
\[
    \HH^i(S^nL) = \HH^i(\rho^*S^nL).
\]

So we have proven the following

\begin{prop}
For a line bundle $L\in\Pic(A)$ one has
\[
    \chi(L_n) = \binom{\chi(L)+n-1}{n}
\]\qed
\end{prop}

\begin{rem}
This result is proven by a somewhat different method in~\cite{lehn99}.
\end{rem}

Next we will attack the cohomology of the restricted bundle
$K^{n-1}L$. The first step in this direction is the following

\begin{lem}
In the notation of diagram~\eqref{isotriviality}, one has
$\nu^*L_n =K^{n-1}L \boxtimes L^n$.
\end{lem}

\begin{proof}

The splitting of the sheaf $\nu^*L_n$ follows from the seesaw
principle (cf.~\cite{mumford70}): For fixed $a\in A$ we have seen
that the restricted morphism $\nu|_{\KA{n-1}\times\{a\}}$ is the
isomorphism wich maps $\KA{n-1}$ to the fiber of the summation
morphism $s$ over the point $na$. Since $\KA{n-1}$ is simply
connected, its Picard group is discrete and it follows that
$\nu^*L_n|_{\KA{n-1}\times\{a\}} \simeq p_{\KA{n-1}}^*
K^{n-1}L|_{\KA{n-1}\times\{a\}}$.

Therefore $\nu^*L_n$ is of the form $K^{n-1}L\boxtimes L_2$ with
$L_2\in \Pic(A)$, and we can compute the component $L_2$ by
considering the restrictions of $\nu$ to $\{\xi_0\}\times A$.

Now complete diagram \eqref{isotriviality} as follows
\begin{equation}\label{singKummerdiagram}
  \begin{CD}
    \KA {n-1} @>{\nu}>> \HA n \\
    @V{\rho'}VV @VV{\rho}V \\
    {K'}^{n - 1} A \times A @>{\nu'}>> \SA n \\
    @V{p_A}VV @VV{\Sigma}V \\
    A @>>n> A.
  \end{CD}
\end{equation}

Here ${K'}^{n-1}A$ denotes the ``singular Kummer'', the fiber over
0 of the addition morphism $\Sigma:\SA{n}\To A$, $\rho' =
\rho|_{{K'}^{n-1}A} \times id_A$ its desingularization, and $\nu'$
is defined analogously to $\nu$.

Now consider a point $\xi_0\in K^{n-1}A$ over $(0,\ldots,0)\in
\SA{n}$. Instead of computing the sheaf
$\nu^*L_n|_{\{\xi_0\}\times A}$ we equivalently compute
${\rho'}^*{\nu'}^*S^nL|_{\{\xi_0\}\times A} = L^n$.
\end{proof}

The next lemma describes the structure of the direct image of
$\shO_A$ under the $n^4$-fold covering $A\stackrel{n}{\To} A$.

\begin{lem}
The direct image $n_*\shO_A$ of the structure sheaf of $A$ splits
into a direct sum of line bundles $L_i, i\in A[n]$, indexed by the
n-torsion points of $A$ and with $L_0=\shO_A$ and
$c_1(L_i)_{\IQ}=0\in \HH^2(A,\IQ)$ for all $i$.
\end{lem}

\begin{proof}
  The splitting of $n_* \shO_A$ is a well known fact,
  cf.~\cite{mumford70}, {\S}7. The triviality of their first rational Chern
  classes follows from the relation:
\[
    L_i^n = L_{in} = L_0 = \shO.
\]

It follows that $nc_1(L_i) = 0$, which proves the lemma.
\end{proof}

Now we have collected all necessary ingredients for the proof of
lemma \ref{RRprop}.

\begin{proof}[Proof of the lemma]
We start again with a line bundle $L$ on $A$ that we twist with a
sufficiently ample bundle $H^N$. By abuse of notation we denote
the resulting bundle by $L$ again. By construction the symmetrized
bundle $S^nL$ is still ample such that the invertible sheaf $L_n$
on $\HA{n}$ as a pull back along the birational
morphism $\rho$ is nef and big.

The same argument shows that the line bundle $K^{n-1}L$ is nef and
big: Using the notations of diagram (\ref{singKummerdiagram}) the
bundle $K^{\prime n-1}L := \nu^{\prime *}S^nL|_{K^{\prime n-1}A}$
is ample and $K^{n-1}L = \rho^{\prime *}K^{\prime n-1}L$ is big
and nef.

So by the Kawamata-Viehweg vanishing theorem (cf.~\cite{viehweg82}),
we have
\[
     \chi(K^{n-1}L) = h^0(K^{n-1}A, K^{n-1}L)
\]

On the one hand, we have due to the K\"{u}nneth formula
\[
    \HH^0(\nu^* L_n) = \HH^0(K^{n-1}L\boxtimes L^n) = \HH^0(K^{n-1}L)
    \otimes \HH^0(L^n).
\]

On the other hand, since $\nu$ is finite and $n$ is a flat
morphism, we can compute $\HH^0(\nu^*L_n)$ alternatively
\[
    \HH^0(\nu^*L_n) = \HH^0(L_n\otimes \nu_*\shO_{K^{n-1}A\times A}) =
    \HH^0(L_n \otimes s^*n_*\shO_A) = \HH^0(L_n\otimes\bigoplus_{i\in A[n]} s^*L_i).
\]

This shows that
\[
    h^0(\KA{n-1}\times A, \nu^*L)  = \sum_{i\in A[n]} h^0(\HA{n}, L_n\otimes
    L_i')
\]

with $L_i' := s^* L_i$.

Since the sheaf $L_n\otimes L_i'$ is still nef and big, the
vanishing theorem of Kawamata and Viehweg implies that
$h^0(\HA{n}, L_n\otimes L_i')$ equals the Euler characteristic of
this line bundle. Therefore, using the classical
Hirzebruch-Riemann-Roch theorem on $\HA{n}$ we have

\begin{align*}
    h^0(\nu^* L_n) &= \sum_{i\in A[n]} \chi(\HA{n},L_n\otimes L_i')\\
                   &= \sum_{i\in A[n]} \int \ch(L_n\otimes L_i')\td(\HA{n})\\
                   &= n^4 \int \ch(L_n)\td(\HA{n}), &&\text{since $c_1(L_i')_\IQ =
                    s^*c_1(L_i)_\IQ =0$}\\
                   &= n^4 \chi(\HA{n}, L_n)\\
                   &= n^4 \dim(S^n\HH^0(L)) &&\text{due to (\ref{H0SnL}).}
\end{align*}

Now combining the computations -- and noting that $h^0(A,L^n) \neq
0$ -- we find

\begin{align*}
    \chi(K^{n-1}A,K^{n-1}L)    &= h^0(K^{n-1}L) =
    \frac{h^0(\nu^* L_n)}{h^0(A,L^n)}
    = \frac{n^4\binom{\frac{c_1(L)^2}{2} + n-1}{n}}{n^2 \frac{c_1(L)^2}{2}}\\
    &= n\binom{\frac{c_1(L)^2}{2} + n-1}{n-1}.
\end{align*}

Once again -- considering the formula as a polynomial in $N$ and
evaluating in $N=0$ -- the formula holds for a general line bundle
$L$.
\end{proof}

\begin{example}
In the case of the Kummer surface $K^1A$, the above formula gives
back the classical Riemann-Roch formula for K3 surfaces:

Remember the diagram
\begin{equation*}
  \begin{CD}
      \KA 1 @<<< \hat A \\
      @V{\epsilon}VV @VVV \\
      A/\mspace{-8mu}\sim @<<p< A.
  \end{CD}
\end{equation*}

Let us start with an symmetric line bundle $L$ on $A$, i.e.\
$L=p^*L'$ with $L'\in \Pic{(A/\mspace{-8mu}\sim)}$. Then $L$
induces a line bundle $M = \eps^*L'$ on $K^1A$. Write
$K^1L=\rho^*\pi_*(L\boxtimes L)^{\mf{S}_2}$, as usual. Then one
has $K^1L=M^2$:

Considering diagram (\ref{ClassCase}), it suffices to show that
$\tilde{\Delta}^*S^2L ={L'}^2$. But this is clear from the
definition of $\tilde{\Delta}$ and $L=p^*L'$.

Now, our HRR formula gives
\[
    \chi(K^1L) = 2\binom{\frac{c_1(L)^2}{2}+1}{1} = c_1(L)^2 + 2
\]

Using that $\eps$ is birational, $p$ is generically 2:1 and the
equality $K^1L=M^2$, one finds

\begin{equation*}
    \chi(K^1L) = \frac{c_1(K^1L)^2}{2} + 2,
\end{equation*}
\vspace{1ex}

which is the classical Riemann-Roch formula for the K3 surface
$K^1A$.
\end{example}


\section{Characteristic numbers of $\KA n$}

In his Ph.D. thesis (cf.~\cite{sawon99}), Sawon calculated all Chern
numbers for the generalized Kummer varieties $\KA n$ with $n \leq 4$.

His method is to use the usual Hirzebruch-Riemann-Roch formula to
express the $\chi_y$-genus, which has been calculated by G\"ottsche
and Soergel (see~\cite{goettsche93}) in terms of characteristic
numbers. This gives him $k$ independent relations between the Chern
numbers on a generalized Kummer variety. These are enough relations
for $n \leq 3$ to determine the values of the Chern numbers. For $n =
4$, one further relation is needed. Sawon managed to calculate the
characteristic numbers $\int_{\KA n} \sqrt{\td(\KA n)}$, which gives
him another relation on the generalized Kummer variety $\KA 4$, thus
enabling him to calculate the Chern numbers of $\KA 4$.

Of course, we can use our formulas to derive relations beetween the
Chern numbers of the generalized Kummer varieties. These relations
together with the relations given by the known value of the
$\chi_y$-genus leads to all Chern numbers of $\KA n$, for $n \leq
5$.

Let $L$ be a line bundle with non-vanishing quadratic form on a
generalized Kummer variety $\KA n$ of complex dimension
$2n$. Let $\lambda := \lambda(L)$ be the characteristic value (as in
section~\ref{s:hrrformula}) of $L$
on $\KA n$. So, $\lambda(L^{\otimes m}) = m^2 \lambda$ for all $m \in
\N_0$.

According to our previous results, we have
\begin{equation}
  \int_{\KA n} \td_{\sqrt{m^2 \lambda/4+1}}(\KA n) = \chi(L^{\otimes m}) =
  (n + 1) \binom{\frac 1 4 m^2 \lambda (n + 1) + n} n
\end{equation}
for $m \in \N_0$. Since both sides are polynomial in $m$, we can
compare the coefficients of $m$, which gives us no more than
$\left\lfloor\frac{n+1} 2\right\rfloor$
independent expressions for some characteristic numbers of $\KA n$
since the coefficients on the left are given by Chern numbers.
For example, the leading coefficients
yield
\begin{equation}
  \int_{\KA n} \sqrt{\td(\KA n)} = \frac{(n + 1)^{n + 1}}{4^n n!},
\end{equation}
a formula which has already been found by Sawon (\cite{sawon99}). The
comparision of coefficients also yields
\begin{equation}
  \label{equ:inttd}
  \int_{\KA n} \td(\KA n) = n + 1,
\end{equation}
which holds on every irreducible symplectic K\"ahler manifold of
complex dimension $2 n$.

We also have to investigate Hirzebruch's
$\chi_y$-genus, which is given by
\begin{equation}
  \label{equ:chiy}
  \chi_y(\KA n) = \chi(\bigwedge_y\nolimits \Omega_X)
  = \sum_{p = 0}^{2n} y^p \int_X \td(X) \ch(\Omega_X^p).
\end{equation}
Here, $\bigwedge_y \Omega_X = \sum_{p = 0}^{2n} y^p \Omega_X^p \in
\KK^0(\KA n)[y]$. Apparently, $\chi_y(\KA n)$ can be expressed by the
Hodge numbers $h^{*, *}(X)$ of $\KA n$, which have been calculated by
G\"ottsche and Soergel (see~\cite{goettsche93}). It is
\begin{equation}
  \chi_y(\KA n) = (n + 1) \sum_{d|(n + 1)} d^3(1 - y + y^2 - \dots
  + (-y)^{\frac{k + 1} d - 1})^2 (-y)^{k + 1 - \frac{k + 1} d}.
\end{equation}
Using~\eqref{equ:chiy}, this yields at most $n$ further
independent expressions for Chern numbers of $\KA n$
(cf.~\cite{sawon99}). One of these expressions is
again~\eqref{equ:inttd}, so summing up, we have at most $[\frac
{n + 1} 2] + n - 1$ linear independent equations for the Chern numbers
on $\KA n$.

We used the computer algebra system Maple to solve the linear
relations for the seven Chern numbers on $\KA 5$ and arrived at the
following table:

\bigskip

\begin{center}
\begin{tabular}{c|r}
  Chern number & Evaluated on $\KA 5$ \\
  \hline
  $c_2^5$ & 84478464 \\
  $c_2^3 c_4$ & 26220672 \\
  $c_2^2 c_6$ & 3141504 \\
  $c_2 c_8$ & 142560 \\
  $c_2 c_4^2$ & 8141472 \\
  $c_4 c_6$ & 979776 \\
  $c_{10}$ & 2592
\end{tabular}
\end{center}

\bigskip

To calculate the ten Chern numbers of $\KA 6$, our methods are not
sufficient since they lead only to eight independent relations.


\appendix

\bibliographystyle{amsplain}
\bibliography{ourbib}

\end{document}